\documentclass[a4paper,11pt]{amsart}
\usepackage{graphicx}
\usepackage{mathptmx}
\usepackage{amsmath}
\usepackage{amssymb}
\usepackage{enumitem}
\usepackage{xcolor}

\newmuskip\pFqmuskip

\newcommand*\pFq[6][8]{%
  \begingroup 
  \pFqmuskip=#1mu\relax
  \mathcode`=\string"8000
  \begingroup\lccode`\~=`\,
  \lowercase{\endgroup\let~}\pFqcomma
  F^{#2}_{#3}{\left(\genfrac..{0pt}{}{#4}{#5}\bigg|#6\right)}%
  \endgroup
}
\newcommand{\pFqcomma}{\mskip\pFqmuskip}

\newtheorem{theorem}{Theorem}

\newtheorem{corollary}[theorem]{Corollary}
\newtheorem{proposition}[theorem]{Proposition}

\begin{document}

\title[A Note on Degenerate Multi-poly-Bernoulli numbers and polynomials]{A Note on Degenerate Multi-poly-Bernoulli numbers and polynomials}

\author{Taekyun  Kim}
\address{Department of Mathematics, Kwangwoon University, Seoul 139-701, Republic of Korea}
\email{tkkim@kw.ac.kr}

\author{DAESAN Kim}
\address{Department of Mathematics, Sogang University, Seoul 121-742, Republic of Korea}
\email{dskim@sogang.ac.kr}


\maketitle

\begin{abstract}
In this paper, we consider the degenerate multi-poly-Bernoulli numbers and polynomials which are defined by means of the multiple polylogarithms and degenerate versions of the multi-poly-Bernoulli numbers and polynomials. We investigate some properties for those numbers and polynomials. In addition, we give some identities and relations for the degenerate multi-poly-Bernoulli numbers and polynomials.
\end{abstract}

\section{Introduction}
For $0 \neq \lambda\in\mathbb{R}$, Carlitz considered the higher-order degenerate Bernoulli polynomials  given by
\begin{equation}
\bigg(\frac{t}{(1+\lambda
t)^{\frac{1}{\lambda}}}\bigg)^{r}(1+\lambda
t)^{\frac{x}{\lambda}}=\sum_{n=0}^{\infty}
\beta_{n,\lambda}^{(r)}(x)\frac{t^{n}}{n!},\quad(\mathrm{see}\
[3]). \label{1}
\end{equation}
When $x=0$, $\beta_{n,\lambda}^{(r)}= \beta_{n,\lambda}^{(r)}(0)$ are called the higher-order degenerate Bernoulli numbers. \par
For $k\in\mathbb{Z}$, the polylogarithm function is defined by
\begin{equation}
    \mathrm{Li}_{k}(x)=\sum_{n=1}^{\infty}\frac{x^{n}}{n^{k}},\quad(\mathrm{see}\ [1,9,17]). \label{2}
\end{equation}
Note that $\displaystyle\mathrm{Li}_{1}(x)=\sum_{n=1}^{\infty}\frac{x^{n}}{n!}=-\log(1-x)\displaystyle$. \\
As is known, the poly-Bernoulli polynomials are defined by
\begin{equation}
\frac{\mathrm{Li}_{k}(1-e^{-t})}{e^{t}-1}e^{xt}=\sum_{n=0}^{\infty}PB_{n}^{(k)}(x)\frac{t^{n}}{n!},\quad(\mathrm{see}\ [9]). \label{3}
\end{equation}
When $x=0$, $PB_n^{(k)}=PB_n^{(k)}(0)$ are called the poly-Bernoulli numbers.
Note that $PB_{n}^{(1)}=B_{n}(x)$, where $B_{n}(x)$ are ordinary Bernoulli polynomials given by
\begin{displaymath}
\sum_{n=0}^{\infty}B_{n}(x)\frac{t^{n}}{n!}=\frac{t}{e^{t}-1}e^{xt},\quad (\mathrm{see}\ [1-21]).
\end{displaymath}
The polyexponential function was introduced by Hardy. Kim-Kim recently introduced a modified version of that, called the modified polyexponential function, which is given by
\begin{equation}
\mathrm{Ei}_{k}(x)=\sum_{n=1}^{\infty}\frac{x^{n}}{(n-1)!n^{k}},\quad(\mathrm{see}\ [8,11,13-15]).\label{4}
\end{equation}
Note that $\mathrm{Ei}_{1}(x)=e^{x}-1$. \par
They also introduced the type 2 poly-Bernoulli polynomials given by
\begin{equation}
\frac{\mathrm{Ei}_{k}\big(\log(1+t)\big)}{e^{t}-1}e^{xt}=\sum_{n=0}^{\infty}B_{n}^{(k)}(x)\frac{t^{n}}{n!},\quad(\mathrm{see}\ [8]). \label{5}
\end{equation}
Note that $B_{n}^{(1)}(x)=B_{n}(x),\ (n\ge 0)$. \par
The degenerate exponential functions are defined by
\begin{equation}
e_{\lambda}^{x}(t)=(1+\lambda t)^{\frac{x}{\lambda}},\quad e_{\lambda}(t)=e_{\lambda}^{1}(t)=(1+\lambda t)^{\frac{1}{\lambda}},\quad(\mathrm{see}\ [8]). \label{6}
\end{equation}
Here we observe that
\begin{equation}
e_{\lambda}^{x}(t)=\sum_{n=0}^{\infty}(x)_{n,\lambda}\frac{t^{n}}{n!},\quad(\mathrm{see}\ [8]),\label{7}
\end{equation}
where $(x)_{0,\lambda}=1$, $(x)_{n,\lambda}=x(x-\lambda)\cdots(x-(n-1)\lambda$, $(n\ge 1)$. \par
It is well known that the Stirling numbers of the second kind are defined by
\begin{equation}
\frac{1}{k!}(e^{t}-1)^{k}=\sum_{n=k}^{\infty}S_{2}(n,k)\frac{t^{n}}{n!},\quad(n \ge 0),\quad(\mathrm{see}\ [8,9,11,14,15]). \label{8}
\end{equation}
For $k_{1},k_{2},\dots,k_{r}\in\mathbb{Z}$, the multiple polylogarithm is defined by
\begin{equation}
\mathrm{Li}_{k_{1},k_{2},\dots,k_{r}}(x)=\sum_{0<n_{1}<n_{2}<\cdots<n_{r}}\frac{x^{n_{r}}}{n_{1}^{k_{1}}n_{2}^{k_{2}}\cdots n_{r}^{k_{r}}},\quad(\mathrm{see}\ [20]), \label{9}
\end{equation}
where the sum is over all integers $n_1,n_2,\dots,n_r$ satisfying $0<n_1<n_2<\cdots<n_r$.

About twenty years ago, the first author introduced the generalized Bernoulli numbers
 $B_{n}^{(k_{1},k_{2},\dots,k_{r})}$ of order $r$ (see [10]) which are given by
\begin{equation}\label{9-1}
\frac{r!\mathrm{Li}_{k_1,k_2,\dots,k_r}(1-e^{-t})}{(e^t -1)^r}=\sum_{n=0}^{\infty}B_n^{(k_1,k_2,\dots,k_r)}\frac{t^n}{n!}.
\end{equation}
Actually, the $r!$ in \eqref{9-1} does not appear in [10]. However, the present definition is
more convenient, since $B_n^{(1,1,\dots,1)}=B_n^{(r)}$ are the Bernoulli numbers of order $r$ (see \eqref{13}).
These numbers would have been called the multi-poly-Bernoulli numbers, since it is a multiple version of
 poly-Bernoulli numbers (see \eqref{3}). Furthermore, we may consider the multi-poly-Bernoulli
 polynomials, which are natural extensions of the multi-poly-Bernoulli numbers, given by
\begin{equation}\label{9-2}
\frac{r!\mathrm{Li}_{k_1,k_2,\dots,k_r}(1-e^{-t})}{(e^t -1)^r}e^{xt}
=\sum_{n=0}^{\infty}B_n^{(k_1,k_2,\dots,k_r)}(x)\frac{t^n}{n!}.
\end{equation}
The multi-poly-Bernoulli polynomials are multiple versions of the poly-Bernoulli polynomials in \eqref{3}.
 We let the interested reader refer to [1] for the detailed properties on those polynomials.\\
\indent In this paper, we consider the degenerate
multi-poly-Bernoulli numbers and polynomials (see \eqref{14})
which are defined by means of the multiple polylogarithms and
degenerate versions of the multi-poly-Bernoulli numbers and
polynomials studied earlier in the literature (see [2]). We
investigate some properties for those numbers and polynomials. In
addition, we give some identities and relations for the degenerate
multi-poly-Bernoulli numbers and polynomials.

\section{Degenerate multi-poly-Bernoulli numbers and polynomials}
From \eqref{9}, we note that
\begin{align}
\frac{d}{dx}\mathrm{Li}_{k_{1},k_{2},\dots,k_{r}}(x)\ &=\ \frac{d}{dx} \sum_{0<n_{1}<n_{2}<\cdots<n_{r}}\frac{x^{n_{r}}}{n_{1}^{k_{1}}n_{2}^{k_{2}}\cdots n_{r}^{k_{r}}}\label{10}  \\
&=\ \frac{1}{x} \sum_{0<n_{1}<n_{2}<\cdots<n_{r}}\frac{x^{n_{r}}}{n_{1}^{k_{1}}n_{2}^{k_{2}}\cdots n_{r}^{k_{r}-1}} \nonumber \\
&=\ \frac{1}{x} \mathrm{Li}_{k_{1},k_{2},\dots,k_{r}-1}(x).\nonumber
\end{align}
Let us take $k_{r}=1$. Then we have
\begin{align}
\frac{d}{dx}\mathrm{Li}_{k_{1},\dots,k_{r-1},1}(x)\ &=\   \sum_{0<n_{1}<n_{2}<\cdots<n_{r-1}}\frac{1}{n_{1}^{k_{1}}n_{2}^{k_{2}}\cdots n_{r-1}^{k_{r-1}}}\sum_{n_{r}=n_{r-1}+1}^{\infty}x^{n_{r}-1}.\label{11}\\
&=\ \frac{1}{1-x}\sum_{0<n_{1}<n_{2}<\cdots<n_{r-1}} \frac{x^{n_{r-1}}}{n_{1}^{k_{1}}n_{2}^{k_{2}}\cdots n_{r-1}^{k_{r-1}}}\nonumber \\
&=\ \frac{1}{1-x} \mathrm{Li}_{k_{1},k_{2},\dots,k_{r-1}}(x) \nonumber
\end{align}
Thus, by \eqref{11}, we get
\begin{equation}
    \mathrm{Li}_{k_{1},1}(x)=\int\frac{1}{1-x} \mathrm{Li}_{k_{1}}(x)dx. \label{12}
\end{equation}
By integration by parts, from \eqref{12}, we note that
\begin{equation*}
    \mathrm{Li}_{1,1}(x)=\int\frac{1}{1-x}(- \log (1-x))dx=\frac{1}{2!}\big(-\log(1-x)\big)^{2}.
\end{equation*}
By induction, we get
\begin{align}
\mathrm{Li}_{\underbrace{1,1,\dots,1}_{r-times}}(x)\ &=\ \frac{(-1)^{r}}{r!}\big(\log(1-x)\big)^{r},\quad (r\in\mathbb{N}), \label{13} \\
&=\ \sum_{l=r}^{\infty}S_{1}(l,r)(-1)^{l-r}\frac{t^{l}}{l!} \nonumber \\
&=\ \sum_{l=r}^{\infty}|S_{1}(l,r)|\frac{t^{l}}{l!} \nonumber
\end{align}
where $S_{1}(l,r)$ (respectively, $|S_{1}(l,r)|$) are the signed (respectively, unsigned) Stirling numbers of the first kind. \par
~~\\
Now, we consider the degenerate multi-poly-Bernoulli polynomials which are degenerate versions of the multi-poly-Bernoulli polynomials in \eqref{9-2} and given by
\begin{equation}
\frac{r!\mathrm{Li}_{k_{1},k_{2},\dots,k_{r}}(1-e^{-t})}{\big(e_{\lambda}(t)-1\big)^{r}}e_{\lambda}^{x}(t)=\sum_{n=0}^{\infty}\beta_{n,\lambda}^{(k_{1},k_{2},\dots,k_{r})}(x)\frac{t^{n}}{n!}. \label{14}
\end{equation}
When $x=0$, $\beta_{n,\lambda}^{(k_{1},k_{2},\dots,k_{r})} =\beta_{n,\lambda}^{(k_{1},k_{2},\dots,k_{r})}(0) $ are called the degenerate multi-poly-Bernoulli numbers. \par
From \eqref{14}, we note that
\begin{align}
\sum_{n=0}^{\infty}\beta_{n,\lambda}^{(k_{1},k_{2},\dots,k_{r})}(x)\frac{t^{n}}{n!}\ &=\ \frac{r!\mathrm{Li}_{k_{1},k_{2},\dots,k_{r}}(1-e^{-t})}{\big(e_{\lambda}(t)-1\big)^{r}}e_{\lambda}^{x}(t) \nonumber \\
&=\ \sum_{l=0}^{\infty}\beta_{l,\lambda}^{(k_{1},k_{2},\dots,k_{r})}\frac{t^{l}}{l!}\sum_{m=0}^{\infty}(x)_{m,\lambda}\frac{t^{m}}{m!} \label{15} \\
&=\ \sum_{n=0}^{\infty}\bigg(\sum_{l=0}^{n}\binom{n}{l}(x)_{n-l,\lambda}\beta_{l,\lambda}^{(k_{1},k_{2},\dots,k_{r})}\bigg)\frac{t^{n}}{n!}. \nonumber
\end{align}
Thus, by \eqref{15}, we get
\begin{equation}
\beta_{n,\lambda}^{(k_{1},k_{2},\dots,k_{r})}(x)=\sum_{l=0}^{n}\binom{n}{l}(x)_{n-l,\lambda} \beta_{l,\lambda}^{(k_{1},k_{2},\dots,k_{r})},\quad (n\ge 0).\label{16}
\end{equation}
From \eqref{1}, \eqref{13} and \eqref{14}, we note that
\begin{equation}
    \beta_{n,\lambda}^{\overbrace{(1,1,\dots,1)}^{r-times}}(x)=\beta_{n,\lambda}^{(r)}(x),\quad(n\ge 0). \label{17}
\end{equation}
\begin{proposition}
For $n\ge 0$, we have
\begin{displaymath}
\beta_{n,\lambda}^{(k_{1},k_{2},\dots,k_{r})}(x)= \sum_{l=0}^{n}\binom{n}{l}(x)_{n-l,\lambda} \beta_{l,\lambda}^{(k_{1},k_{2},\dots,k_{r})},
\end{displaymath}
\begin{displaymath}
\beta_{n,\lambda}^{\overbrace{(1,1,\dots,1)}^{r-times}}(x)=\beta_{n,\lambda}^{(r)}(x).
\end{displaymath}
\end{proposition}
From \eqref{14}, we note that
\begin{align}
&\sum_{n=0}^{\infty}\beta_{n,\lambda}^{(k_{1},\dots,k_{r})}(x)\frac{t^{n}}{n!}=\frac{r!}{\big(e_{\lambda}(t)-1\big)^{r}}\mathrm{Li}_{k_{1},\dots,k_{r}}\big(1-e^{-t}\big)e_{\lambda}^{x}(t) \label{18} \\
&=\ \frac{r!}{\big(e_{\lambda}(t)-1\big)^{r}}\sum_{0<n_{1}<\cdots<n_{r-1}}\frac{1}{n_{1}^{k_{1}}\cdots n_{r-1}^{k_{r-1}}}\sum_{n_{r}=n_{r-1}+1}^{\infty}\frac{\big(1-e^{-t}\big)^{n_{r}}}{n_{r}^{k_{r}}}e_{\lambda}^{x}(t) \nonumber \\
&=\ \frac{r!}{\big(e_{\lambda}(t)-1\big)^{r}}\sum_{0<n_{1}<\cdots<n_{r-1}}\frac{\big(1-e^{-t}\big)^{n_{r-1}}}{n_{1}^{k_{1}}\cdots n_{r-1}^{k_{r-1}}}\sum_{n_{r}=1}^{\infty}\frac{\big(1-e^{-t}\big)^{n_{r}}}{(n_{r}+n_{r-1})^{k_{r}}}e_{\lambda}^{x}(t) \nonumber \\
&=\ \frac{r!e_{\lambda}^{x}(t)}{\big(e_{\lambda}(t)-1\big)^{r}}\sum_{0<n_{1}<\cdots<n_{r-1}}\frac{\big(1-e^{-t}\big)^{n_{r-1}}}{n_{1}^{k_{1}}\cdots n_{r-1}^{k_{r-1}}}\sum_{n_{r}=1}^{\infty}\frac{n_{r}!}{(n_{r}+n_{r-1})^{k_{r}}}\sum_{l=n_{r}}^{\infty}(-1)^{l-n_r}S_{2}(l,n_{r})\frac{t^{l}}{l!} \nonumber
\end{align}
To proceed further, we observe that, for any integer $k$, we have
\begin{equation}\label{19}
(x+y)^{-k}=\sum_{m=0}^{\infty}(-1)^m \binom{k+m-1}{m} x^{-k-m} y^m.
\end{equation}
Now, from \eqref{18} and \eqref{19}, we have

\begin{align*}
\sum_{n=0}^{\infty}\beta_{n,\lambda}^{(k_{1},\dots,k_{r})}(x)\frac{t^{n}}{n!}
&=\ \frac{rt}{e_{\lambda}(t)-1}\bigg(\frac{(r-1)!e_{\lambda}^{x}(t)}{\big(e_{\lambda}(t)-1\big)^{r-1}}\sum_{0<n_{1}<\cdots<n_{r-1}} \frac{\big(1-e^{-t}\big)^{n_{r-1}}}{n_{1}^{k_{1}}\cdots n_{r-1}^{k_{r-1}-m}}\bigg)\nonumber \\
&\quad \times\frac{1}{t}\sum_{l=1}^{\infty}\bigg(\sum_{n_{r}=1}^{l}n_{r}!(-1)^{l-n_{r}} S_{2}(l,n_{r})\sum_{m=0}^{\infty}\binom{k_{r}+m-1}{m}(-1)^{m}n_{r}^{-k_{r}-m}\bigg)\frac{t^{l}}{l!}\nonumber\\
&=\ \frac{rt}{e_{\lambda}(t)-1}\sum_{m=0}^{\infty}\binom{k_{r}+m-1}{m}(-1)^{m}\sum_{j=0}^{\infty}\beta_{j,\lambda}^{(k_{1},\dots,k_{r-1}-m)}(x)\frac{t^{j}}{j!}\nonumber \\
&\quad \times \sum_{l=0}^{\infty}\bigg(\sum_{n_{r}=1}^{l+1}\frac{n_{r}!(-1)^{l-n_{r}-1} S_{2}(l+1,n_{r})}{l+1}n_{r}^{-k_{r}-m}\bigg)\frac{t^{l}}{l!}\nonumber\\
&=\ \frac{rt}{e_{\lambda}(t)-1}\sum_{k=0}^{\infty}\bigg(\sum_{l=0}^{k}\sum_{n_{r}=1}^{l+1}\sum_{m=0}^{\infty}\binom{k_{r}+m-1}{m}(-1)^{m}\binom{k}{l}\nonumber\\
&\quad\times\frac{n_{r}!(-1)^{l-n_{r}-1} S_{2}(l+1,n_{r})}{l+1} n_{r}^{-k_{r}-m}\beta_{k-l,\lambda}^{(k_{1},k_{2},\dots,k_{r-1}-m)}(x)\bigg) \frac{t^{k}}{k!}\nonumber
\end{align*}
\begin{align}\label{20}
&=\ r\sum_{p=0}^{\infty}\beta_{p,\lambda}\frac{t^{p}}{p!}\sum_{k=0}^{\infty}\bigg( \sum_{l=0}^{k}\sum_{n_{r}=1}^{l+1}\sum_{m=0}^{\infty}\binom{k_{r}+m-1}{m}(-1)^{m}\nonumber\\
&\quad\times\binom{k}{l}\frac{n_{r}!(-1)^{l-n_{r}-1} S_{2}(l+1,n_{r})}{l+1} n_{r}^{-k_{r}-m}\beta_{k-l,\lambda}^{(k_{1},k_{2},\dots,k_{r-1}-m)}\bigg)\frac{t^{k}}{k!}\nonumber \\
&=\ \sum_{n=0}^{\infty}\bigg(\sum_{k=0}^{n}\sum_{l=0}^{k}\sum_{n_{r}=1}^{l+1}\sum_{m=0}^{\infty}
r\binom{k_{r}+m-1}{m}\binom{k}{l}\binom{n}{k}(-1)^{m}\nonumber \\
&\quad\times\frac{n_{r}!(-1)^{l-n_{r}-1} S_{2}(l+1,n_{r})}{l+1} n_{r}^{-k_{r}-m}\beta_{k-l,\lambda}^{(k_{1},k_{2},\dots,k_{r-1}-m)}(x)\beta_{n-k,\lambda}\bigg)\frac{t^{n}}{n!},
\end{align}
where $\beta_{n,\lambda}$ are the Carlitz's degenerate Bernoulli numbers with $\beta_{n,\lambda}^{(1)}=\beta_{n,\lambda}$. Therefore, by \eqref{20}, we obtain the following theorem.

\begin{theorem}
For $k_{1},k_{2},\dots,k_{r}\in\mathbb{Z}$ and $n\ge 0$, we have
\begin{align}\label{21}
\beta_{n,\lambda}^{(k_{1},k_{2},\dots,k_{r})}(x)\ &=\ r \sum_{k=0}^{n}\sum_{l=0}^{k}\sum_{n_{r}=1}^{l+1}\sum_{m=0}^{\infty}
\binom{k_{r}+m-1}{m}\binom{k}{l}\binom{n}{k}(-1)^{m}\nonumber \\
&\quad\times\frac{n_{r}!(-1)^{l-n_{r}-1} S_{2}(l+1,n_{r})}{l+1} n_{r}^{-k_{r}-m}\beta_{n-k,\lambda}\beta_{k-l,\lambda}^{(k_{1},k_{2},\dots,k_{r-1}-m)}(x).
\end{align}
\end{theorem}

Replacing $k_r$ by $-k_r$ in \eqref{21}, we obtain the following corollary.
\begin{corollary}
For $k_{1},k_{2},\dots,k_{r}\in\mathbb{Z}$ and $n\ge 0$, we have
\begin{align*}
\beta_{n,\lambda}^{(k_{1},k_{2},\dots,-k_{r})}(x)\ &=\ r \sum_{k=0}^{n}\sum_{l=0}^{k}\sum_{n_{r}=1}^{l+1}\sum_{m=0}^{\infty}
\binom{k_{r}}{m}\binom{k}{l}\binom{n}{k}\nonumber \\
&\quad\times\frac{n_{r}!(-1)^{l-n_{r}-1} S_{2}(l+1,n_{r})}{l+1} n_{r}^{k_{r}-m}\beta_{n-k,\lambda}\beta_{k-l,\lambda}^{(k_{1},k_{2},\dots,k_{r-1}-m)}(x).
\end{align*}
\end{corollary}

From \eqref{14}, we have
\begin{displaymath}
    \sum_{n=0}^{\infty}\bigg(\beta_{n,\lambda}^{(k_{1},\dots,k_{r})}(x+1)-\beta_{n,\lambda}^{(k_{1},\dots,k_{r})}(x)\bigg)\frac{t^{n}}{n!}=\frac{r!e_{\lambda}^{x}(t)}{\big(e_{\lambda}(t)-1\big)^{r-1}}\mathrm{Li}_{k_{1},\dots,k_{r}}\big(1-e^{-t}\big).
\end{displaymath}
\begin{align*}
    &=\frac{r!e_{\lambda}^{x}(t)}{\big(e_{\lambda}(t)-1\big)^{r-1}}\sum_{0<n_{1}<\cdots<n_{r-1}}\frac{1}{n_{1}^{k_{1}}\cdots n_{r-1}^{k_{r-1}}}\sum_{n_r=n_{r-1}+1}^{\infty}\frac{\big(1-e^{-t}\big)^{n_{r}}}{n_{r}^{k_{r}}}\\
    &=\frac{r!}{\big(e_{\lambda}(t)-1\big)^{r-1}} e_{\lambda}^{x}(t)\sum_{0<n_{1}<\cdots<n_{r-1}}\frac{\big(1-e^{-t}\big)^{n_{r-1}}}{n_{1}^{k_{1}}\cdots n_{r-1}^{k_{r-1}}}\sum_{n_r=1}^{\infty}\frac{\big(1-e^{-t}\big)^{n_{r}}}{(n_{r}+n_{r-1})^{k_{r}}}\\
    &=r\sum_{m=0}^{\infty}\binom{k_{r}+m-1}{m}(-1)^{m}\frac{(r-1)!e_{\lambda}^{x}(t)}{\big(e_{\lambda}(t)-1\big)^{r-1}}\mathrm{Li}_{k_{1},\dots,k_{r-1}-m}\big(1-e^{-t}\big)\\
    &\qquad\quad\times\sum_{n_{r}=1}^{\infty}n_{r}^{-k_{r}-m}n_{r}!\sum_{l=n_{r}}^{\infty}(-1)^{l-n_r}S_{2}(l,n_{r})\frac{t^{l}}{l!}\\
    &= r\sum_{m=0}^{\infty}\binom{k_{r}+m-1}{m}(-1)^{m}\sum_{j=0}^{\infty}\beta_{j,\lambda}^{(k_{1},\dots,k_{r-2},k_{r-1}-m)}(x)\frac{t^{j}}{j!}\\
    &\qquad\quad\times\sum_{l=1}^{\infty}\bigg(\sum_{n_r=1}^{l}n_{r}^{-k_{r}-m}n_{r}!(-1)^{l-n_r}S_{2}(l,n_{r})\bigg)\frac{t^{l}}{l!}
\end{align*}
\begin{align*}
    &=\ \sum_{n=1}^{\infty}\bigg(r\sum_{m=0}^{\infty}\binom{k_{r}+m-1}{m}(-1)^{m}\sum_{l=1}^{n}\sum_{n_{r}=1}^{l}\binom{n}{l}n_{r}^{-k_{r}-m}n_{r}! \\
&\quad\times (-1)^{l-n_r}S_{2}(l,n_{r})\beta_{n-l,\lambda}^{(k_{1},k_{2},\dots,k_{r-2},k_{r-1}-m)}(x)\bigg)\frac{t^{n}}{n!}.
\end{align*}
Therefore, we obtain the following theorem.
\begin{theorem}
For $n,r\ge 1$ and $k_{1},k_{2},\dots,k_{r}\in\mathbb{Z}$, we have
\begin{align*}
\frac{1}{r}\bigg(\beta_{n,\lambda}^{(k_{1},\dots,k_{r})}(x+1)-\beta_{n,\lambda}^{(k_{1},\dots,k_{r})}(x)\bigg)&=\sum_{m=0}^{\infty}\binom{k_{r}+m-1}{m}(-1)^{m}\sum_{l=1}^{n}\sum_{n_{r}=1}^{l}\binom{n}{l}n_{r}^{-k_{r}-m}\\
&\quad\times n_{r}!(-1)^{l-n_r}S_{2}(l,n_{r})\beta^{(k_{1},\dots,k_{r-2},k_{r-1}-m)}_{n-l,\lambda}(x).
\end{align*}
\end{theorem}
By the definition of degenerate multi-poly-Bernoulli polynomials, we get
\begin{displaymath}
    \sum_{n=0}^{\infty}\beta_{n,\lambda}^{(k_{1},k_{2},\dots,k_{r})}(x+y)\frac{t^{n}}{n!}=\frac{r!}{\big(e_{\lambda}(t)-1\big)^{r}}\mathrm{Li}_{k_{1},k_{2},\dots,k_{r}}(1-e^{-t})e_{\lambda}^{x+y}(t).
\end{displaymath}
\begin{align*}
    &=\ \frac{r!}{\big(e_{\lambda}(t)-1\big)^{r}}\mathrm{Li}_{k_{1},k_{2},\dots,k_{r}}(1-e^{-t})e_{\lambda}^{x}(t)e_{\lambda}^{y}(t) \\
    &=\ \sum_{l=0}^{\infty}\beta_{l,\lambda}^{(k_{1},k_{2},\dots,k_{r})}(x)\frac{t^{l}}{l!}\sum_{m=0}^{\infty}(y)_{m,\lambda}\frac{t^{m}}{m!} \\
    &=\ \sum_{n=0}^{\infty}\bigg(\sum_{l=0}^{n}\binom{n}{l}\beta_{l,\lambda}^{(k_{1},\dots,k_{r})}(x)(y)_{n-l,\lambda}\bigg)\frac{t^{n}}{n!}.
\end{align*}
Thus, we note that
\begin{displaymath}
    \beta_{n,\lambda}^{(k_{1},\dots,k_{r})}(x+y)=\sum_{l=0}^{n}\binom{n}{l}\beta_{l,\lambda}^{(k_{1},\dots,k_{r})}(x)(y)_{n-l,\lambda},\quad(n\ge 0).
\end{displaymath}
Finally, we note that $\beta_{n,\lambda}^{(k_{1},\dots,k_{r})}(x)$ is not a Sheffer sequence.

\section{Conclusion}
In [3], Carlitz initiated study of degenerate versions of Bernoulli and Euler polynomials, namely the degenerate Bernoulli and Euler polynomials. In recent years, some mathematicians intensively studied various versions of many special numbers and polynomials and quite a few interesting results were found about them (see [6,7,9,11-16,18,20]). Here we would like to mention that study of degenerate versions is not limited only to polynomials but can be extended also to transcendental functions like gamma functions (see [12]). \\
\indent In this paper, we considered the degenerate multi-poly-Bernoulli numbers and polynomials which are defined by means of the multiple polylogarithms. They are degenerate versions of the multi-poly-Bernoulli numbers and polynomials, and multiple versions of the degenerate poly-Bernoulli numbers and polynomials ($r=1$ case of \eqref{14}). We investigated some properties for those numbers and polynomials. In fact, among other things, we derived some explicit expressions of the degenerate multi-poly-Bernoulli numbers and polynomials.\\
\indent As it turns out, exploration for degenerate versions of some special polynomials and transcendental functions has been very rewarding and fruitful. Indeed, the study of degenerate versions has applications to differential equations, identities of symmetry and probability theory as well as to number theory and combinatorics. We let the reader refer to the Conclusion section of [14] for the details. It is one of our future projects to continue this line of research and find applications not only in mathematics but also in science and engineering.

\end{document}